\documentclass[11pt,reqno]{article}
\usepackage{amsmath, amssymb}

\newtheorem{theorem}{Theorem}[section]

\newtheorem{corollary}[theorem]{Corollary}
\newtheorem{lemma}[theorem]{Lemma}

\def \proof {\noindent {\bf Proof.}\ \ }
\def \remark {\noindent {\bf Remark.}\ \ }

\def \endproof {{\mbox{}\nolinebreak\hfill\rule{2mm}{2mm}\par\medbreak}}

\def \R {\mathbb{R}}

\def \E {\mathbb{E}}
\def \P {\mathbb{P}}

\def \a {\alpha}

\def \e {\varepsilon}

\def \d {\delta}

\def \l {\lambda}

\def \s {\sigma}

\def \< {\langle}
\def \> {\rangle}
\def \^ {\widehat}

\def \Prob {{\rm Prob}}
\def \id {{\it id}}

\begin{document}
\title {Frame expansions with erasures: an approach through 
        the non-commutative operator theory}
\author{Roman Vershynin\footnote{
Partially supported by the 
New Faculty Research Grant of the University of California-Davis 
and by the NSF grant DMS 0401032}}

\maketitle

\begin{abstract}
In modern communication systems such as the Internet, random losses of information 
can be mitigated by oversampling the source. This is equivalent to expanding 
the source using overcomplete systems of vectors (frames), as opposed to 
the traditional basis expansions. 
Dependencies among the coefficients in frame expansions 
often allow for better performance comparing to bases 
under random losses of coefficients. 
We show that for any $n$-dimensional frame, 
any source can be linearly reconstructed from only $\sim n \log n$ randomly chosen 
frame coefficients, with a small error and with high probability.
Thus every frame expansion withstands random losses better (for worst
case sources) than the orthogonal basis expansion, 
for which the $n \log n$ bound is attained.
The proof reduces to M.Rudelson's selection theorem on random 
vectors in the isotropic position, which is based on the non-commutative
Khinchine's inequality. 
\end{abstract}

\section{Introduction}

Representation of signals using frames, which are overcomplete sets of vectors, 
is advantageus over basis expansions in a variety of practical applications.  
Dependencies among the coefficients of the overcomplete representations guarantee 
a better stability in presence of noise, quantization, erasures,
as well as greater freedom of design comparing to bases. This general paradigm 
is confirmed by many experiments and some theoretical work, see e.g.
\cite{D}, \cite{G2}, \cite{GVT}, \cite{GKK}, \cite{KDG}, \cite{BO}, \cite{CK}
and the bibliography contained therein. 

Of particular importance are the dependencies contained in frame expansions for 
design of communication systems. The redundancy of frames can mitigate random 
losses of expansion coefficients that occur in packet-based 
communication systems such as the Internet. Detection and retransmission  
of lost packets in such systems takes much longer than their 
successful transmission. This is main source of delays known to all network users. 
Such delays are unacceptable for many applications, such as the real-time video.
It is thus esirable for the receiver to be able to approximately reconstruct 
the information sent to him from {\em whatever} packets he receives, despite 
the loss of some packets. There should exist certain dependencies among the 
packets, otherwise the information contained in a missing packet would be 
irrevocably lost. Then, what is the best way to distribute the information 
among the packets so that each packet is equally important? 
Equivalently, this is the problem of the Multiple Description Coding (MDC) theory, 
where one wishes to communicate information over a set of 
parallel channels, each of which either works perfectly or not at all.

The idea originated in \cite{GKK} was to use frame expansions
to distribute the information among the packets with some dependencies.
One can view this communication scheme as follows:
\begin{equation}						\label{comm scheme}
x \in \R^n
\xrightarrow{}     \boxed{\substack{\text{frame} \\ \text{expansion}}}
\xrightarrow{y \in \R^m}   \boxed{\substack{\text{transmission} \\ \text{(losses)}}}
\xrightarrow{\^ y \in \R^k}   \boxed{\text{{\footnotesize reconstruction}}}
\xrightarrow{}
\^ x \in \R^n
\end{equation}
The source information is viewed as a vector $x \in \R^n$.
This vector is represented by its $m \ge n$ expansion coefficients with respect to 
some fixed frame. These coefficients are sent over the network in $m$ packets, 
each in its own packet. 
Due to unpredictable communication losses, the user 
receives only a random subset of these packets, say $k$ in average. 
The user applies the linear reconstruction to the received coefficients in hope 
that the reconstruction error would be small with graceful probability.
The fundamental problem is\footnote{In this paper, we neglect the 
quantization issues, which are treated in \cite{GVT} and \cite{GKK})}: 

\begin{quote}
{\em How many random coefficients of a frame expansion does the user 
need to receive to be able to linearly reconstruct the source vector 
with a small error and with large probability?}
\end{quote}

The work on this question, both theoretical and experimental, was initiated 
in \cite{GKK} and continued in \cite{KDG} and \cite{CK},
see also a survey paper \cite{G2}.
Both cases were considered: $k < n$, which clearly requires a statistical 
model of input vector $x$, and $k \ge n$. The performance of 
the frame representations was compared to that of the classical block 
channel-coded basis representations.

In the present paper we look for a best bound on $k$ 
which works for {\em all} frames and {\em all} source vectors $x$. 
Does every frame necessarily perform better than the trivial 
frame, the orthonormal basis -- or, more generally, an orthonormal 
basis in $\R^n$ each of whose elements is repeated $s$ times? 
Communicating a source vector $x$ with the trivial frame is equivalent to 
sending each of the $n$ coefficients of the orthonormal expansion of $x$ precisely 
$s$ times. To be able to reconstruct $x$, 
the user must receive each of the $n$ coefficients at least once. 
This is possible with probability at least $1-\e$ 
only if the user receives $k \ge C(\e) n \log n$ random coefficients in total. 
This gives the lower bound on $k$ in the question above. Remarkably, the upper
bound matches.

\begin{theorem}							\label{main}
  For any uniform tight frame in $\R^n$ and any source vector $x$, 
  the linear reconstruction from $k$ random coefficients of $x$ yields an 
  approximation error at most $\e$ with probability $1-\e$, provided 
  $k \ge C(\e) n \log n$.
\end{theorem}
Here $C(\e)$ is a constant that depends only on $\e$; this dependence is
discussed in the next section.
Tightness of a frame is assumed here only for simplicity. 

Note that the optimal bound on $k$ does not depend on the size $m$ of the frame, 
so there may be many lost coefficients -- in fact, most of them may be lost.
Hence it is not the number of the lost coefficients that determines 
the performance but the number $k$ of received coefficients.

As argued in \cite{G2}, one advantage of frame representations over 
the traditional block channel-coded basis representations is 
that frames allow for a real time reconstruction of the source.
The receiver can attempt to reconstruct a source vector -- such as a still image
or video -- in real time as the packets arrive, starting from the very first 
successfully received coefficient.
Within one communication session, the number of received coefficients $k$ will 
thus grow in time from $1$ to possibly $m$, and the quality of reconstruction will 
improve as more coefficients arrive. 
(In contrast to this, in the block channel-coded bases model the user must wait 
until $n$ coefficients arrive). 
Theorem \ref{main} states that, with {\em any} frame design and {\em any} source, 
the reconstruction quality will reach a nearly optimal level 
as soon as $\sim n \log n$ coefficients are received, 
so one may stop the session then.

Theorem \ref{main} shows that every frame must withstand random losses better 
than the trivial frame, the one formed by repeating the elements
of the orhogonal basis. Of course, there exist frames that perform 
better than the trivial frame. The problem of optimal desing of such frames
is addressed in \cite{GKK} and \cite{CK}. 
As noticed e.g. in \cite{GVT}, a set of $m=sn$ random points $(x_i)$
taken independently with the uniform distribution on the unit sphere 
$S^{n-1}$ forms a frame which approaches a tight frame with large probability, 
provided the redundancy $s \to \infty$. 
Consequently, a random $k$-element subset of this set also 
forms an almost tight frame with large probability, provided $k \ge tn$ and 
$t$ is large. Then one can linearly reconstruct any source vector $x$
from using its $k$ random coefficients with respect to the frame $(x_i)$
with probability $1-\e$, provided $k \ge C(\e) n$. Hence for this frame, 
the logarithmic factor is not needed in the number received coefficients $k$.

Our proof of Theorem \ref{main} is based on a result of M.Rudelson in the
asymptotic convex geometry about vectors in the isotropic position \cite{R2}. 
There exists a remarkable equivalence of the theories.
All of the following classes coincide in $\R^n$ (up to an appropriate rescaling), 
see \cite{V}:
\begin{itemize}
  \item the class of tight frames,
  \item the class of contact points of convex bodies,
  \item the class of John's decompositions of the identity,
  \item the class of vectors in the isotropic position.
\end{itemize}  
The selection theorem of M.Rudelson \cite{R2} can thus be interpreted 
as a result about frames, which leads to Theorem \ref{main}.
In order to obtain an exponentially large probability in Theorem \ref{main}
and because of a slightly different model of random selection in 
M.Rudelson's theorem, we will prove the latter with some necessary modifications.
Two proofs of Rudelson's theorem are known. The one which was
historically the first \cite{R1} uses majorizing measures, a deep technique in modern 
probability theory developed by M.Talagrand (see \cite{T}). 
The other proof \cite{R2} is the one we follow in the present paper.
It is based on the noncommutative operator theory, 
more precisely on the noncommutative Khinchine's inequality due to F.Lust-Piquard 
and G.Pisier (see \cite{LP}, \cite{P}, \cite{R2}).

Section \ref{s:frames} relates the frames to the decompositions of the 
identity and offers a precise form of Theorem \ref{main}. 
Section \ref{s:noncommutative} discusses the noncommutative Khinchine's
inequality and Pisier's proof of Rudelson's lemma.
In Section \ref{s:proof} we show how Rudelson's lemma implies
a precise form of Theorem \ref{main}.

\section{Frames as decompositions of identity and their random parts}     \label{s:frames}

For an introduction to frames, see \cite{D} and \cite{C}.
A system of vectors $(x_i)$ finite or infitite, in a Hilbert space,  
is called a {\em frame} if there exist $A > 0$ and $B > 0$ 
(the {\em frame bounds}) such that 
$$
A \|x\|^2  \le  \sum_i |\< x, x_i\> |^2  \le  B \|x\|^2
\ \ \ \text{holds for all $x \in \R^n$.}
$$

Our Hilbert space will be $\R^n$ with its canonical scalar product.
We will specialize to {\em uniform frames}, those for which $\|x_i\| = 1$
for all $i$, and to {\em tight frames}, for which $A = B$.
The reason for considering only tight frames is the simple fact
that a frame has frame bounds $(A,B)$ if and only if it is 
$\sqrt{AB}$-equivalent to some tight frame (see \cite{C}). 
By being $M$-equivalent we mean 
that there exists a linear operator $T$ that maps elements of one frame to 
the other with $\|T\|\|T^{-1}\| \le M$.

We will view frame expansions as decompositions of identity. 
A pair of vectors $(x,y)$ in $\R^n$ defines a one-dimensional 
linear operator $x \otimes y$ given by $(x \otimes y)(z) = \< x,z\> y$. 
Then for any system of vectors $(x_i)_{i=1}^m$ with $\|x_i\| = 1$
and for the identity operator $\id$ on $\R^n$ one has
\begin{equation}						\label{frame vs decomp}
\text{$(x_i)_{i=1}^m$ is a uniform tight frame in $\R^n$ if and only if} \ \ \
\id = \frac{n}{m} \sum_{i=1}^m x_i \otimes x_i.
\end{equation}
Communication scheme \eqref{comm scheme} based on a uniform tight frame 
$(x_i)_{i=1}^m$ works as follows. A source vector $x \in \R^n$ is represented
through the expansion \eqref{frame vs decomp}, i.e. 
$$
x = \frac{n}{m} \sum_{i=1}^m \< x_i, x \> x_i,
$$
and the coefficients $y(i) := \< x_i, x \> $, $i=1,\ldots,m$ are sent over 
the network. At each given time during the communication session, the user
has received a random subset $\s \subset\{1,\ldots,m\}$ of these coefficients.
The user applies to them the linear reconstruction, computing
\begin{equation}						\label{xhat}
\^ x = \frac{n}{|\s|} \sum_{i \in \s} \< x_i, x \> x_i
\end{equation}
in hope that the error $\| x - \^ x\|$ would be small with large probability.
The question is -- how large should $|\s|$ for this to hold?

More formally, the random subset $\s$ is realized by including each element of 
$\{1,\ldots,m\}$ into $\s$ independently with probability $k/m$, 
where $0 < k < m$ is some fixed number. Then $\s$ is a {\em random subset 
of $\{1,\ldots,m\}$ of average size $k$}.

\begin{theorem}							\label{main full}
  Let $(x_i)_{i=1}^m$ be a uniform tight frame in $\R^n$, and $\e > 0$.
  Let $\s$ be a random subset of $\{1,\ldots,m\}$ of average size 
  $k \ge C (n/\e^2) \log (n/\e^2)$. Then 
  $$
  \P \Big\{ \Big\| \id - \frac{n}{|\s|} \sum_{i \in \s} x_i \otimes x_i \Big\| 
            > \e t \Big\}
  \le C e^{-t^2}
  $$
  in the (only interesting) range $0 < t < 1/\e$.
\end{theorem}
Here and thereafter $C, C_1, \ldots$ denote absolute constants, 
whose values for convenience may be different from line to line
(but they do not depend on anything). 

Theorem \ref{main full} gives an asymptotically optimal bound on the 
required number $k$ of received coefficients 
in communication scheme \eqref{comm scheme}:

\begin{corollary}						\label{main cor}
  Let $(x_i)_{i=1}^m$ be a uniform tight frame in $\R^n$. 
  Let $\e \in (0,1)$, $t > 1$ and $k \ge C (n/\e^2) \log (n/\e^2)$.
  With probability at least $1 - Ce^{-t^2}$, the linear 
  reconstruction \eqref{xhat} from a random subset $\s$ 
  of average size $k$ gives the error
  $$
  \| x - \^ x\| < \e t \ \ \ 
  \text{for all possible sources $x \in \R^n$.}
  $$
\end{corollary}
\endproof

\noindent Theorem \ref{main} clearly follows from Corollary \ref{main cor}.

\qquad

\remark The proof also shows that the average approximation error 
in Theorem~\ref{main cor} is small, 
$\E \| x - \^ x\| < \e$.

\section{Noncommutative Khinchine's inequality and Rudelson's theorem}  
							\label{s:noncommutative}

The main ingredient in the proof of Theorem \ref{main full}
is the following result of M.Rudelson \cite{R2}.

\begin{lemma}[M.Rudelson]					\label{mark}
  Let $(z_i)$ be a finite collection of vectors in $\R^d$. Then  
  $$
  \Big( \E \Big\| \sum_i \e_i z_i \otimes z_i \Big\|^p \Big)^{1/p}
  \le C (p + \log d)^{1/2} \max_i \|z_i\|
      \cdot \Big\| \sum_i z_i \otimes z_i \Big\|^{1/2}.
  $$
\end{lemma}

G.Pisier (\cite{P}, see \cite{R2}) discovered an approach to this result 
via the noncommutative operator theory, which greatly simplified the original 
proof of M.Rudelson \cite{R1}. For completeness, we give a proof of Lemma \ref{mark}
since only case $p=1$ was treated explicitely in the literature.

Lemma \ref{mark} reduces to the noncommutative Khinchine inequality 
due to F.Lust-Piquard and G.Pisier (see \cite{LP}, \cite{P}, \cite{R2}).
In the noncommutative operator theory, the role of scalars is 
played by linear operators. Beside the usual operator norm, an
operator $Z$ on $\R^d$ has the norm in the Schatten class $C_p^d$ for $p \ge 1$,
defined as follows. Let $s_i(Z)$ be the $s$-numbers of $Z$, that is
the eigenvalues of $Z^*Z$. The norm in the Schatten class is then 
$\|Z\|_{C_p^d} = (\sum_{i=1}^d s_i(Z)^p)^{1/p}$.

\begin{theorem}[Non-commutative Khinchine's inequality 
               \cite{LP}, \cite{P}, see \cite{R2}]        \label{noncom}
  Let $2 \le p < \infty$.
  For any finite sequence $(Z_i)$ in $C_p^d$ one has
  $$
  R((Z_i)) 
    \le \Big( \E \Big\| \sum_i \e_i Z_i \Big\|_{C_p^d}^p \Big)^{1/p}    
    \le  C \sqrt{p} \cdot R((Z_i)),
  $$
  where 
  $$
  R((Z_i))  =  \max \Big( 
    \Big\| \big( \sum Z_i^* Z_i \big)^{1/2} \Big\|_{C_p^d}, \ 
    \Big\| \big( \sum Z_i Z_i^* \big)^{1/2} \Big\|_{C_p^d}
  \Big).
  $$
\end{theorem}

In the scalar case, that is for $d=1$, Teorem \ref{noncom} is the 
classical Khinchine's inequality (see e.g. \cite{LT} Lemma 4.1).

\noindent {\bf Proof of Lemma \ref{mark}. }
Note that for every $r \ge 1$ and every operator $Z \in C_r^d$,
$$
\|Z\|_{C_r^d} = \Big( \sum_{i=1}^d s_i(Z)^r \Big)^{1/r}
\le d^{1/r} \max_i s_i(Z).
$$
Let $r = p + \log d$. Then $d^{1/r} \le e$, hence 
\begin{equation}						\label{norms equiv}
\|Z\| \le \|Z\|_{C_r^d} \le e \|Z\|.
\end{equation}
We apply the noncommutative Khinchine's inequality for 
$Z_i = z_i \otimes z_i$. Note that 
$Z_i^* Z_i = Z_i Z_i^* = \|z_i\|^2 z_i \otimes z_i$.
By \eqref{norms equiv},
\begin{align*}
\Big( \E \Big\| \sum_i \e_i z_i \otimes z_i \Big\|^p \Big)^{1/p}
&\le \Big( \E \Big\| \sum_i \e_i z_i \otimes z_i \Big\|_{C_r^d}^p \Big)^{1/p} \\
&\le C \sqrt{r} \Big\| \Big( \sum_i \|z_i\|^2 z_i \otimes z_i \Big)^{1/2}
             \Big\|_{C_r^d}^p \\
&\le C e \sqrt{r} \max_i \|z_i\| \cdot \Big\| \Big( \sum_i z_i \otimes z_i 
            \Big)^{1/2} \Big\|.
\end{align*}
In view of our choice of $r$, this completes the proof of Lemma \ref{mark}.
\endproof

\section{Proof of Theorem~\ref{main full}} 		\label{s:proof}

\paragraph{Moments and tails}
The tail probability in Theorem \ref{main full}
can be computed by estimating the moments. This is described in the 
following standard lemma.
For an $\a \ge 1$, the $\psi_\a$-norm of a random variable $Z$
is defined as 
$$
\|Z\|_{\psi_\a} 
  = \inf \big\{ \l > 0 : \; \E \exp|Z/\l|^\a  \le e \big\}.
$$  

\begin{lemma} (see \cite{LT} Lemmae 3.7 and 4.10)		\label{moments tails}
  Let $Z$ be a nonnegative random variable, and let $\a = d/2$
  for some positive integer $d$.
  The following are equivalent: 

  (i) there exists a constant $K>0$ such that 
  $$
  (\E Z^p)^{1/p} \le K p^\a
  \ \ \ \text{for all $p \ge 2$};
  $$
  
  (ii) there exists a constant $K>0$ such that 
  $$
  \P \{ Z > Kt \} \le 2 \exp (-t^{1/\a})
  \ \ \ \text{for all $t > 0$};
  $$

  (iii) there exists a constant $K>0$ such that 
  $$
  \|Z\|_{\psi_\a} \le K.
  $$
Furthermore, the constants in (i), (ii) and (iii) depend only on $\a$
and on each other.
\end{lemma}

\begin{corollary}          				 \label{pnorm}
  Let $Z$ is a nonegative random variable and $p \ge 2$. Then
  $$
  (\E Z^p)^{1/p}  \le  C p \log ( \E \exp Z )
  $$
  for all $p \ge 1$.
\end{corollary}

\proof
Let $M = \|Z\|_{\psi_1}$. 
Assume first that $M \ge 1$. We have
$$
\E \exp (Z / M) = e.
$$
By Lemma \ref{moments tails}, $(\E (Z/M)^p)^{1/p}  \le  C p$.
Then by Jensen's inequality
\begin{eqnarray*}
(\E Z^p)^{1/p}
  &\le&  C p M 
    =    C p M \log ( \E \exp (Z /M) ) \\
  & = &  C p \log \big( \E \exp (Z /M) \big)^M \\
  &\le&  C p \log ( \E \exp Z).
\end{eqnarray*}

For a general nonnegative variable $Z$, note that 
$\|1+Z\|_{\psi_1} \ge 1$, hence by the previous argument
$$
(\E Z^p)^{1/p} 
\le (\E (1+Z)^p)^{1/p}
\le C p \log ( \E \exp (1+Z) )
= Cep \log ( \E \exp Z).
$$
This completes the proof.
\endproof

\paragraph{Symmetrization}
We start our proof of Theorem \ref{main full} 
with the decomposition \eqref{frame vs decomp},
$$
x = \frac{n}{m} \sum_{i=1}^m \< x_i, x \> x_i.
$$
To realize a random subset $\s$, we introduce selectors $(\d_i)_{i=1}^m$,
that is independent $\{0,1\}$-valued random variables with means
$\E \d_i = \d$, where $\d = \frac{k}{m}$.
Then $\s = \{i :\; \d_1 = 1\}$ is a random subset of $\{1,\ldots,m\}$ of 
average size $k$. 

Disregarding for a moment a difference between the random size $|\s|$ 
and its mean $k$, thanks to Lemma \ref{moments tails}
we can compute the probability estimate in Theorem \ref{main full}
by estimating the moments 
$$
E_p = \Big( \E \Big\| \id - \frac{n}{k} \sum_{i \in \s} x_i \otimes x_i 
               \Big\|^p \Big)^{1/p}
= \Big( \E \Big\| \id - \frac{n}{k} \sum_{i=1}^m \d_i x_i \otimes x_i 
               \Big\|^p \Big)^{1/p}
$$ 
for $p \ge 2$.
This will be done in several steps. 

At the first step, we apply the classical symmetrization tecnique 
(see \cite{LT} 6.2). We look at 
$Y = \id - \frac{n}{k} \sum_{i=1}^m \d_i x_i \otimes x_i$
as a random variable (random operator) and consider its independent copy $Y'$.
Since $\E Y' = 0$, Jensen's inequality yields $\E\|Y\|^p \le \E\|Y-Y'\|^p$, 
hence 
$$
E_p \le \Big( \E \Big\| \frac{n}{k} \sum_{i=1}^m (\d_i - \d'_i) x_i \otimes x_i 
               \Big\|^p \Big)^{1/p},
$$ 
where $(\d'_i)_{i=1}^m$ is an independent copy of $(\d_i)_{i=1}^m$.
Let $(\e_i)$ be a sequence of independent symmetric $\{-1,1\}$-valued random variables, 
independent of both $(\d_i)$ and $(\d'_i)$. 
Since $\d_i - \d'_i$ is a symmetric random variable, it is distributed
identically to $\e_i(\d_i - \d'_i)$. By Minkowski's inequality, 
\begin{align}							\label{ep sym}
E_p 
&\le \Big( \E \Big\| \Big( \frac{n}{k} \sum_{i=1}^m \e_i \d_i x_i \otimes x_i \Big)
                   - \Big( \frac{n}{k} \sum_{i=1}^m \e_i \d'_i x_i \otimes x_i \Big)
               \Big\|^p \Big)^{1/p} \nonumber\\
&\le 2 \Big( \E \Big\| \frac{n}{k} \sum_{i=1}^m \e_i \d_i x_i \otimes x_i 
               \Big\|^p \Big)^{1/p}.
\end{align}

\paragraph{Bounding the moments}
Let us fix a realization of the selectors $(d_i)$ (hence a set $\s$)
and denote by $\E_\e$ the expectation with respect to $(\e_i)$.
The number of nonzero elements among $z_i = \d_i x_i$, $i = 1,\ldots, m$
is $d = |\s| = \sum_{i=1}^m \d_i$. Consequently, we can view $z_i$
as vectors in $\R^d$. Applying Lemma \ref{mark} to them, we obtain
\begin{align*}
\Big( \E_\e \Big\| \frac{n}{k} \sum_{i=1}^m \e_i \d_i x_i \otimes x_i 
               \Big\|^p \Big)^{1/p}
&= \frac{n}{k} \Big( \E_\e \Big\| \sum_{i \in \s} \e_i z_i \otimes z_i 
               \Big\|^p \Big)^{1/p} \\
&\le \frac{Cn}{k} (p + \log|\s|)^{1/2} 
    \Big\| \sum_{i=1}^m \d_i x_i \otimes x_i \Big\|^{1/2}
\end{align*}
By \eqref{ep sym} and Cauchy-Schwartz inequality, we get
\begin{align}								\label{Ep split}
E_p 
&\le 2 \Big( \E \, \E_\e \Big\| \frac{n}{k} \sum_{i=1}^m \e_i \d_i x_i \otimes x_i 
               \Big\|^p \Big)^{1/p} \nonumber\\
&\le 2C \sqrt{\frac{n}{k}}
    \big[ \E (p+\log|\s|)^p \big]^{1/2p}
    \left[ \E \Big\| \frac{n}{k} \sum_{i=1}^m \d_i x_i \otimes x_i \Big\|^p
          \right]^{1/2p}.
\end{align}
The first expectation in \eqref{Ep split} is estimated by 
Minkowski's inequality and Corollary~\ref{pnorm} as
\begin{align*}
\big[ \E (p + \log|\s|)^p \big]^{1/2p}
&\le \big[ p + (\E \log^p |\s|)^{1/p} \big]^{1/2} \\
&\le \big[ p + C p \log \E|\s| \big]^{1/2} 
= [p + C p \log k]^{1/2}
\le C (p \log k)^{1/2}.
\end{align*}
The second expectation in \eqref{Ep split} is estimated
by Minkowski's inequality as
$$
\left[ \E \Big\| \frac{n}{k} \sum_{i=1}^m \d_i x_i \otimes x_i \Big\|^p
          \right]^{1/2p}
\le (1 + E_p)^{1/2}.
$$
Summarizing, \eqref{Ep split} becomes
$$
E_p^2 \le C p \Big( \frac{n \log k}{k} \Big) (1+E_p).
$$
Denoting $a = \frac{n \log k}{k}$ and solving for $E_p$, we have
$$
E_p \le C(ap + \sqrt{ap}),
$$
thus
$$
\min(E_p, 1) \le C \sqrt{ap}.
$$
Since $E_p = (\E Z^p)^{1/p}$ for 
$Z = \| \id - \frac{n}{k} \sum_{i \in \s} x_i \otimes x_i \|$,
we have 
$$
\big[ \E (\min(Z,1))^p \big]^{1/p}
\le \min (E_p, 1)
\le C \sqrt{ap}.
$$
By Corollary \ref{moments tails}, 
\begin{equation}						\label{minZ1}
\P \{ \min(Z,1) > C_1 \sqrt{a} \, t \}
\le 2 \exp(-t^2)
\ \ \ \text{for all $t > 0$.}
\end{equation}
Now recall the restriction on $k$ in Theorem \ref{main full}, 
$k \ge C (n/\e^2) \log (n/\e^2)$. 
By choosing $C$ large enough, we can make
$$
C_1 \sqrt{a} = C_1 \sqrt{\frac{n \log k}{k}}
\le \e/10.
$$
In view of the definition of $Z$, \eqref{minZ1} implies
\begin{equation}						\label{almost}
\P \Big\{ \Big\| \id - \frac{n}{k} \sum_{i \in \s} x_i \otimes x_i 
          \Big\| 
         > \frac{\e t}{10} \Big\}
\le 2 \exp(-t^2)
\ \ \ \text{for all $0 < t < 10/\e$.}
\end{equation}

\paragraph{Replacing the average size of the random set by its actual size}
It remains to replace $k$ by $|\s|$ in \eqref{almost}.
Indeed, since $|\s| = \sum_{i=1}^m \d_i$ is a sum of $m$ independent 
$\{ 0, 1\}$-valued random variables $\d_j$ with $\E \d_j = \d = \frac{k}{m}$,
Bernstein's inequality (see \cite{Pe}) shows that 
for $s \le 2 \d m = 2k$ one has
$$
\Prob \big\{  \big| |\s| - k \big| > s \big\}
\le  2\exp \Big( - \frac{s^2}{8 \d m}
           \Big)   
\le  2\exp \Big( - \frac{s^2}{8 k} \Big).
$$
Then for $s = \frac{\e t k}{10}$, 
$$
\Prob \Big\{  \Big| \frac{|\s|}{k} - 1 \Big| > \frac{\e t}{10} \Big\}
\le  2\exp \Big(- \frac{\e^2 t^2 k}{800} \Big)    
\le  2\exp(-t^2).     
$$
If both events $| \frac{|\s|}{k} - 1 | \le \frac{\e t}{10}$
and $\| \id - \frac{n}{k} \sum_{i \in \s} x_i \otimes x_i \| 
\le \frac{\e t}{10}$ hold, which happens with probability 
at least $1 - 4 \exp(-t^2)$, then by the triangle inequality 
$\| \frac{n}{k} \sum_{i \in \s} x_i \otimes x_i \| 
\le 1 + \frac{\e t}{10} < 2$, hence 
$$
\Big\| \id - \frac{n}{|\s|} \sum_{i \in \s} x_i \otimes x_i \Big\| 
\le \Big\| \id - \frac{n}{k} \sum_{i \in \s} x_i \otimes x_i \Big\| 
  + \Big\| \Big(1 - \frac{k}{|\s|} \Big) 
           \frac{n}{k} \sum_{i \in \s} x_i \otimes x_i \Big\| 
\le \frac{\e t}{10} + \frac{4 \e t}{10}
< \e t.
$$
Thus $k$ may be replaced by $|\s|$ in \eqref{almost} 
at the cost of replacing $\frac{\e t}{10}$ by $\e t$.
This completes the proof of Theorem \ref{main full}.
\endproof

\qquad

\noindent Department of Mathematics \\
University of California \\
Davis, CA 95616 \\
U.S.A. \\
E-mail: {\small\tt%
vershynin@math.ucdavis.edu} 

\end{document}